\documentclass[a4paper,twoside,10pt]{article}
\usepackage{amsmath}
\usepackage{amsthm}
\usepackage{amsfonts}
\usepackage{amssymb, amscd}
\usepackage[dvips]{graphicx}
\usepackage[all]{xy}
\usepackage{eucal}
\numberwithin{equation}{section}
\title{Quasi-Deformations of $\sll$ using twisted derivations}
\author{Daniel Larsson and
    Sergei D. Silvestrov\\
    \footnotesize \emph{Centre for Mathematical Sciences, Department of
   Mathematics, Lund Institute}\\
   \footnotesize \emph{of Technology, Lund University, Box
   118, SE-221 00 Lund, Sweden}\\
   \footnotesize dlarsson@maths.lth.se, ssilvest@maths.lth.se}
\theoremstyle{definition}
\newtheorem{dfn}{Definition}
\newtheorem{example}{Example}

\theoremstyle{plain}
\newtheorem{remark}{Remark}
\theoremstyle{plain}
\newtheorem{prop}{Proposition}
\newtheorem{thm}[prop]{Theorem}

\def\F{\mathbb{F} }

\def\C{\mathbb{C} }

\def\Z{\mathbb{Z} }
\def\N{\mathbb{N} }
\def\mfD{\mathfrak{D} }
\def\g{\mathfrak{g} }

\def\A{\mathcal A }

\def\mfD{ \mathfrak{D} }

\def\ps{ \partial_\sigma }
\def\sll{\mathfrak{sl}_2(\F)}
\def\Jacksll{{}_q\mathfrak{sl}_2(\F)}

\def\A{\mathcal A }

\def\cs{\circlearrowleft }

 \DeclareMathOperator{\id}{id}

\DeclareMathOperator{\Ann}{Ann}

\DeclareMathOperator{\Lin}{\mathcal{L}}\DeclareMathOperator{\linSpan}{LinSpan}

\begin{document}
\maketitle
\begin{abstract}
In this paper we apply a method devised in
\cite{HartLarsSilv1D,LarsSilv1D} to the three-dimensional simple
Lie algebra $\sll$. One of the main points of this deformation
method is that the deformed algebra comes endowed with a canonical
twisted Jacobi identity. We show in the present paper that when
our deformation scheme is applied to $\sll$ we can, by choosing
parameters suitably, deform $\sll$ into the Heisenberg Lie algebra
and some other three-dimensional Lie algebras in addition to more
exotic types of algebras, this being in stark contrast to the
classical deformation schemes where $\sll$ is rigid. The resulting
algebras are quadratic and we point out possible connections to
``geometric quadratic algebras'' such as the Artin--Schelter
regular algebras, studied extensively since the beginning of the
90's in connection with non-commutative projective geometry.
\end{abstract}
\section{Introduction}
The idea to deform algebraic, analytic and geometric structures
within the appropriate category is obviously not new. The first
modern appearance is often attributed to Kodaira and Spencer
\cite{KodairaSpencerD} and deformations of complex structures on
complex manifolds. This was, however, soon extended and
generalized in an algebraic-homological setting by Gerstenhaber,
Grothendieck and Schlessinger. But the idea to deform mathematical
structures and objects certainly traces back even further: the
Taylor polynomial of a (holomorphic) function can be viewed as a
natural deformation of the function to a finite polynomial
expression. Nowadays deformation-theoretic ideas permeate most
aspects of both mathematics and physics and cut to the very core
of theoretical and computational problems. In the case of Lie
algebras, which will be our primary concern, quantum deformations
(or $q$-deformations) and quantum groups associated to Lie
algebras have been in style for over twenty years, still growing
richer by the minute. This area began a period of rapid expansion
around 1985 when Drinfel'd \cite{Drin1D} and Jimbo \cite{Jimbo1D}
independently considered deformations of $\mathbf{U}(\g)$, the
universal enveloping algebra of a Lie algebra $\g$, motivated,
among other things, by their applications to the Yang--Baxter
equation and quantum inverse scattering methods \cite{KulResD}.
Since then several other versions of ($q$-) deformed Lie algebras
have appeared, especially in physical contexts such as string
theory. The main objects for these deformations were
infinite-dimensional algebras, primarily the Heisenberg algebras
(oscillator algebras) and the Virasoro algebra, see
\cite{ChaiKuLukD,CurtrZachos1D,DamKuD,HartLarsSilv1D,HelSil-bookD}
and the references therein. For more details how these algebras
are important in physics, see
\cite{ChaiKuLukD,CurtrZachos1D,DamKuD,QuantFieldStringsD,DiFranMathiSenCD,Fuchs1D,Fuchs2D},
for instance. We note that the deformed objects in the above cases
seldom, if ever, belong to the original category of Lie algebras.
However, one retains the undeformed objects in the appropriate
limit. The deformations we will consider do not necessarily have
this important property as we will see. Therefore we refer to
these as \emph{quasi-deformations} thereby emphasizing this
crucial difference explicitly. Strictly speaking we do not even
consider deformations in the classical sense of
Gerstenhaber--Grothendieck--Schlessinger. Instead, suppose
$\mathfrak{g}$ is the Lie algebra we want to deform and suppose
that $\mathfrak{g}\xrightarrow{\rho}\mathfrak{gl}(\A)$ is a
representation of $\mathfrak{g}$ on a commutative, associative
algebra $\A$ with unity, where $\mathfrak{gl}(\A)$ is the Lie
algebra (under the commutator bracket) of linear operators on the
underlying vector space of $\A$. Then our deformation scheme can
be diagrammatically depicted as
$$
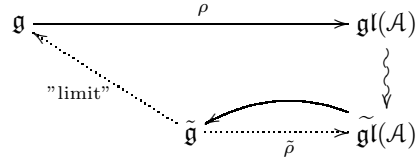
\xymatrix{\mathfrak{g}\ar[rrrr]^{\rho}&&&&\mathfrak{gl}(\A)\ar@{~>}[d]\\
    &&\tilde{\mathfrak{g}}\ar@{.>}[ull]^{\text{"limit"}}\ar@{.>}[rr]_{\tilde\rho}&&
    \widetilde{\mathfrak{gl}}(\A)\ar@/_1pc/[ll]}$$
where the "squiggly" line is the "deformation" procedure of
substituting the original operators with deformed
($\sigma$-twisted, with $\sigma$ an algebra endomorphism on $\A$)
versions. The dotted arrows indicate the fact that we do not
necessarily come back to the object we started with when going
back the way the arrays point, i.e., the very reason for calling
it ``quasi-deformation''. The algebra $\tilde{\mathfrak{g}}$ is
then to be considered as the "deformation" of $\mathfrak{g}$. So
what we actually change, or deform, is the given representation of
$\mathfrak{g}$ and the bracket product in $\mathfrak{gl}(\A)$.
This explains why the tilde $\,\widetilde{}\,$ appears above
$\mathfrak{gl}(\A)$ in the bottom row of the above diagram. We
note that the above method could maybe be generalized to other
algebras besides Lie algebras. However, we do not yet know how, or
to what extent this is possible as of this moment.

The basic undeformed object ($\mathfrak{g}$ from the above
diagram) in this article is the classical $\mathfrak{sl}_2(\C)$,
the Lie algebra generated as a vector space by elements $H,E$ and
$F$ with relations
\begin{align}\label{eq:sl2}
    \langle H,E\rangle = 2E, \qquad\langle H,F\rangle = -2F,
    \qquad\langle E,F\rangle = H.
\end{align} This Lie algebra is simple and perhaps the single most
important one since any (complex) semi-simple Lie algebra includes
a number of copies of $\mathfrak{sl}_2(\C)$ (see for instance, \cite{Serre1D},
page 43--44). It can thus be argued that this is also the most
important algebra to deform. In addition, $\mathfrak{sl}_2(\C)$ has a
lot of interesting representations. One such, which is our basic starting
point, is the following in terms of first order differential
operators acting on some vector space of functions in the variable
$t$:
\begin{align*}
    E &\mapsto \partial,\qquad H\mapsto -2t\partial,\qquad F \mapsto
    -t^2\partial.
\end{align*}
This is what we will generalize to first order operators acting on
an algebra $\A$, where $\partial$ is replaced by $\ps$, a
$\sigma$-derivation on $\A$ (see Section \ref{sec:qhlsigma} for
definitions). In \cite{HartLarsSilv1D} it was shown that, for
a $\sigma$-derivation $\ps$ on a commutative associative
algebra $\A$ with unity, the rank one module $\A\cdot\ps$
admits the structure of a $\C$-algebra with a $\sigma$-deformed
commutator
$$\langle a\cdot\ps,b\cdot\ps\rangle=\sigma(a)\cdot\ps(b\cdot\ps)-
\sigma(b)\cdot\ps(a\cdot\ps)$$ as a multiplication, satisfying a
generalized twisted six-term Jacobi identity.

The present article builds on this and provides an elaborated
example of a quasi-deformed $\sll$ (where $\mathbb{F}$ is a field
of characteristic zero) by the above outlined method. The result
becomes a natural quasi-deformation of $\sll$ and we will show
that it is a qhl-algebra in general (see Section
\ref{sec:qhlsigma} for the definition). By choosing different base
algebras $\A$ we obviously get different algebra structures on
$\A\cdot\ps$. So in a way we have two different deformation
"parameters", namely $\sigma$ and $\A$ (of course, changing $\A$
and not changing $\sigma$ is mathematically absurd, since $\sigma$
is dependent on $\A$, but it is nice to loosely think about $\A$
as an independent deformation parameter).

Since the defining relations for $\mathfrak{sl}_2(\C)$ are
quadratic, passing to the universal enveloping algebra gives us a
quadratic algebra with non-homogenous relations. The deformations
of $\mathfrak{sl}_2(\C)$ we consider also yield quadratic algebras
but in general the defining relations are far more involved.
However, in some cases we obtain algebras which, if not already
explicitly studied in the literature, then strongly resembling
such, for instance the Sklyanin algebra and more generally, the
Artin--Schelter regular algebras. For this reason it is natural to
suspect that some non-commutative geometry, such as point- and
line-modules, could be involved even in some of our algebras.

The paper is organized as follows: in Section
\ref{sec:qhlsigma} we recall the relevant definitions and results
from the papers \cite{HartLarsSilv1D} and \cite{LarsSilv1D}; in
Section \ref{sec:defsl2} we deform $\sll$ with $\A$ as general as
can be in this situation and we deduce some necessary conditions
for everything to make sense; in Section \ref{sec:k[t]} we take
the algebra $\A$ to be simply $\F[t]$ and calculate some
properties of the corresponding deformations. Finally, in Section
\ref{sec:k[t]/t^3} we consider the deformations arising when
$\A=\F[t]/(t^N)$ for $N=3$ and for general non-negative integer
$N$, showing that we get a family of algebras parameterized by the
non-negative integers generating deformations at roots of
unity.
\section{Qhl-algebras associated with
$\sigma$-derivations}\label{sec:qhlsigma} We now fix notation and
state the main definitions and results from \cite{HartLarsSilv1D}
and \cite{LarsSilv1D} needed in this paper.

Throughout we let $\F$ denote a field of characteristic zero and
 $\A$ be a commutative, associative $\F$-algebra with unity
$1$. Furthermore $\sigma$ will denote an endomorphism on $\A$.
Then by a \emph{twisted derivation} or \emph{$\sigma$-derivation}
on $\A$ we mean an $\F$-linear map $\ps:\A\to\A$ such that a
$\sigma$-twisted Leibniz rule holds:
\begin{align}\label{eq:sigmaLeibniz}
\ps(ab)=\ps(a)b+\sigma(a)\ps(b).
\end{align}
Among the best known $\sigma$-derivations are:
\begin{itemize}
    \item[$\bullet$] $(\partial\, a)(t)=a'(t)$, \emph{the ordinary differential
    operator} with the ordinary Leibniz rule, i.e., $\sigma=\id$.
    \item[$\bullet$] $(\ps\, a)(t)=a(t+1)-a(t)$, \emph{the
        shifted difference operator}; \\
    $\sigma$-Leibniz: $(\ps\, (ab))(t)=(\ps
    a)(t)b(t)+a(t+1)(\ps b)(t).$
    In this case $\sigma=\mathbf{s}$, where $\mathbf{s}(f)(t):=f(t+1)$.
    \item[$\bullet$] $(\ps\, a)(t)=(D_q a)(t)$, \emph{the Jackson
        $q$-derivation operator}; \\
    $\sigma$-Leibniz: $(D_q\, (ab))(t)=(D_q
    a)(t)b(t)+a(qt)(D_q b)(t).$
    Here $\sigma=\mathbf{t}_q$, where $\mathbf{t}_q f(t):=f(qt).$
\end{itemize}In the paper \cite{HartLarsSilv1D}
the notion of a hom-Lie algebra as a deformed version of a Lie
algebra was introduced, motivated by some of the examples of
deformations of the Witt and Virasoro algebras constructed using
$\sigma$-derivations.
\begin{dfn}\label{def:homLieD}
    Let $L$ be a vector space. A \emph{hom-Lie algebra} structure
    on $L$ is a linear map $\alpha:L\to L$ and a bracket
    $\langle\cdot,\cdot\rangle$ such that
    \begin{itemize}
    \item[$\bullet$] $\langle x,y\rangle =-\langle y,x\rangle$
    \item[$\bullet$] $\cs_{x,y,z}\langle (\alpha+\id)(x),\langle
        y,z\rangle\rangle=0,$
    \end{itemize} for $x,y,z\in L$, where $\cs_{x,y,z}$ denotes
    cyclic summation with respect to $x,y,z$. A morphism of
    hom-Lie algebras $(L,\alpha)\xrightarrow{\phi}(L',\alpha')$ is
    an algebra homomorphism satisfying the intertwining condition
    $\phi\circ\alpha=\alpha'\circ\phi$.
\end{dfn} Note that for $\alpha=\id$ we retain the definition of
a Lie algebra. However, finding examples of more general kinds of
deformations associated to $\sigma$-derivations, prompted the
introduction in \cite{LarsSilv1D} of the following structure
generalizing hom-Lie algebras.
\begin{dfn}\label{def:qhl}
    Let $L$ be a vector space and let $\alpha,\beta$ be linear
    maps on $L$. Also let $\omega_L:D_\omega\to\Lin_\F(L)$ be a map
    from $D_\omega\subseteq L\times
    L$ to $\Lin_\F(L)$ the space of linear maps on $L$ over $\F$. Then $L$ is a
    \emph{quasi-hom-Lie algebra} or a \emph{qhl-algebra} equipped with
    a bilinear product
    $\langle\cdot,\cdot\rangle : L\times L\to L$ if
    the following conditions hold:
    \begin{itemize}
        \item[$\bullet$] $\langle
        \alpha(x),\alpha(y)\rangle=\beta\circ\alpha\langle
        x,y\rangle$
        \item[$\bullet$] $\langle x,y\rangle =\omega_L(x,y)\langle y,x\rangle$
        \item[$\bullet$] $\cs_{x,y,z} \omega_L(z,x)(\langle \alpha(x),\langle
       y,z\rangle\rangle+\beta\langle x,\langle
       y,z\rangle\rangle)=0,$
    \end{itemize}for $x,y,z\in L$ and $(z,x), (x,y), (y,z)\in D_\omega\subseteq L\times L$.
\end{dfn}Note that with this definition quasi-hom-Lie algebras include not only hom-Lie algebras as a
subclass, but also color Lie algebras and in particular Lie
super-algebras \cite{LarsSilv1D}. We let $\mfD_\sigma(\A)$ denote the
set of $\sigma$-derivations on $\A$. Fixing a homomorphism
$\sigma:\A\to\A$, an element $\ps\in\mfD_\sigma(\A)$, and an element
$\delta\in\A$, we assume that these objects satisfy the following two
conditions:
\begin{align}\label{eq:GenWittCond1}
    &\sigma(\Ann(\ps))\subseteq \Ann(\ps),\\
    &\ps(\sigma(a)) = \delta\sigma(\ps(a)),\quad\text{for }a\in\A,\label{eq:GenWittCond2}
\end{align} where $\Ann(\ps):=\{a\in\A\,|\, a\cdot
\ps=0\}$. Let $\A\cdot\ps=\{a\cdot\ps\;|\;a\in\A\}$ denote the
cyclic $\A$-submodule of $\mfD_\sigma(\A)$ generated by $\ps$ and
extend $\sigma$ to $\A\cdot\ps$ by
$\sigma(a\cdot\ps)=\sigma(a)\cdot\ps$. The following theorem, from
\cite{HartLarsSilv1D}, introducing an $\F$-algebra structure on
$\A\cdot\ps$ making it a quasi-hom-Lie algebra, is of central
importance for the present paper.
\begin{thm}\label{thm:GenWitt}
If {\rm (\ref{eq:GenWittCond1})} holds then the map
$\langle\cdot,\cdot\rangle_\sigma$ defined by setting
\begin{align} \label{eq:GenWittProdDef}
    \langle a\cdot\ps,b\cdot\ps\rangle_\sigma=(\sigma(a)\cdot\ps)\circ(b\cdot\ps)-(\sigma(b)\cdot\ps)
    \circ(a\cdot\ps),
\end{align}for $a,b\in\A$ and where $\circ$ denotes composition of maps, is a well-defined
$\F$-algebra product on the $\F$-linear space $\A\cdot\ps$. It
satisfies the following identities for $a,b,c\in\A${\rm :}
\begin{align}\label{eq:GenWittProdFormula}
    &\langle a\cdot\ps,
    b\cdot\ps\rangle_\sigma=(\sigma(a)\ps(b)-\sigma(b)\ps(a))\cdot\ps,\\
    &\langle a\cdot\ps, b\cdot\ps\rangle_\sigma=-\langle b\cdot\ps,
    a\cdot\ps\rangle_\sigma,\label{eq:GenWittSkew}
\end{align}
and if, in addition, {\rm (\ref{eq:GenWittCond2})} holds,
we have the deformed six-term Jacobi identity
\begin{equation} \label{eq:GenWittJacobi}
    \cs_{a,b,c}\,\big(\langle\sigma(a)\cdot\ps,\langle
    b\cdot\ps,c\cdot\ps
\rangle_\sigma\rangle_\sigma+\delta\cdot\langle a\cdot\ps, \langle
b\cdot\ps,c\cdot\ps\rangle_\sigma\rangle_\sigma\big)=0.
\end{equation}
\end{thm}The algebra $\A\cdot\ps$ from the theorem is then a
qhl-algebra with $\alpha=\sigma$, $\beta=\delta$ and
$\omega=-\id_{\A\cdot\ps}$.
For the reader's convenience we sketch the proof of this result. For
full details see \cite{HartLarsSilv1D}.
\begin{proof}
First of all, skew-symmetry and bilinearity are obvious from the
definition (\ref{eq:GenWittProdDef}). The bracket is well-defined
by the assumption (\ref{eq:GenWittCond1}) on $\Ann(\ps)$. Closure
of $\langle\cdot,\cdot\rangle_\sigma$ follows from
\begin{multline*}
    \langle a\cdot\ps, b\cdot\ps\rangle_\sigma(c)=
    (\sigma(a)\cdot\ps)((b\cdot\ps)(c))-(\sigma(b)\cdot\ps)((a\cdot\ps)(c))=\\
    =(\sigma(a)\ps(b)-\sigma(b)\ps(a))\cdot\ps(c)+(\sigma(a)\sigma(b)-\sigma(b)\sigma(a))
    \cdot\ps(\ps(c))
\end{multline*}since $(\sigma(a)\sigma(b)-\sigma(b)\sigma(a))\cdot\ps(\ps(c))=0$ by the commutativity of $\A$.
Let us now prove the deformed Jacobi identity
(\ref{eq:GenWittJacobi}).
This follows from the commutativity of $\A$, the assumption (\ref{eq:GenWittCond2})
and the
fact that $\ps$ is a $\sigma$-derivation,
by the following cyclic summations:
\begin{multline}\label{eq:GenWittCalc1}
    \cs_{a,b,c}\langle \sigma(a)\cdot\ps,\langle
    b\cdot\ps,c\cdot\ps\rangle_\sigma\rangle_\sigma=\\
   =\cs_{a,b,c}\langle\sigma(a)\cdot\ps,(\sigma(b)\ps(c)-\sigma(c)\ps(b))\cdot\ps
   \rangle_\sigma=\\
    =\cs_{a,b,c}\big
    (\sigma^2(a)\sigma^2(b)\ps^2(c)-\sigma^2(a)\sigma^2(c)
    \ps^2(b)\big )\cdot\ps+\\
     +\cs_{a,b,c}\big
    (\sigma^2(c)\sigma(\ps(b))\ps^2(\sigma(a))-\sigma^2(b)
    \sigma(\ps(c))\ps(\sigma(a))\big)\cdot\ps+\\
    +\cs_{a,b,c}\big
    (\sigma^2(a)\ps(\sigma(b))\ps(c)-\sigma^2(a)\ps(\sigma(c))
    \ps(b)\big )\cdot\ps=\\
    =\cs_{a,b,c}\big
    (\sigma^2(a)\ps(\sigma(b))\ps(c)-\sigma^2(a)\ps(\sigma(c))
    \ps(b)\big )\cdot\ps.
\end{multline}The first two terms vanished when adding up cyclically.
Re-write the equality (\ref{eq:GenWittProdFormula}) as $\langle
b\cdot\ps,c\cdot\ps\rangle_\sigma=(\ps(c)\sigma(b)-\ps(b)\sigma(c))\cdot\ps$
which is possible since $\A$ is commutative. Now the second part
of (\ref{eq:GenWittJacobi}) becomes with the aid of
(\ref{eq:GenWittCond2}):
\begin{multline*}
    \cs_{a,b,c}\delta\cdot\langle a\cdot\ps,\langle
    b\cdot\ps,c\cdot\ps\rangle_\sigma\rangle_\sigma=\\
    =\cs_{a,b,c}\delta\cdot\langle
   a\cdot\ps,(\ps(c)\sigma(b)-\ps(b)\sigma(c))\cdot\ps\rangle_\sigma=\\
    =\cs_{a,b,c}\big
      (\delta\cdot\sigma(a)\ps^2(c)\sigma(b)-
      \delta\cdot\sigma(a)\ps^2(b)\sigma(c)\big )\cdot\ps+\\
      +\cs_{a,b,c}\big (-\ps(\sigma(c))\sigma^2(b)\ps(a)+
      \ps(\sigma(b)\sigma^2(c)\ps(a) \big)\cdot\ps=\\
      =\cs_{a,b,c}\big (-\ps(\sigma(c))\sigma^2(b)\ps(a)+
      \ps(\sigma(b)\sigma^2(c)\ps(a) \big)\cdot\ps,
\end{multline*}where the first term becomes zero because of cyclic
summation and the commutativity of $\A$. Combining this with
(\ref{eq:GenWittCalc1}) yields (\ref{eq:GenWittJacobi}).
\end{proof}
\section{Quasi-Deformations}\label{sec:defsl2} Let $\A$ be a commutative,
associative $\F$-algebra with unity $1$, $t$ an element of $\A$,
and let $\sigma$ denote an $\F$-algebra endomorphism on $\A$.
Also, let $\mfD_\sigma(\A)$ denote the linear space of
$\sigma$-derivations on $\A$. Choose an element $\ps$ of
$\mfD_\sigma(\A)$ and consider the $\F$-subspace $\A\cdot\ps$ of
elements on the form $a\cdot\ps$ for $a\in\A$. We will usually
denote $a\cdot\ps$ simply by $a\ps$. Notice that $\A\cdot\ps$ is a
left $\A$-module. By Theorem \ref{thm:GenWitt} there is a
skew-symmetric algebra structure on this $\A$-module given by
\begin{multline}\label{eq:bracket}\langle
a\cdot\ps, b\cdot\ps\rangle=\sigma(a)\cdot\ps(b\cdot\ps)-
\sigma(b)\cdot\ps(a\cdot\ps)=\\
=\sigma(a\ps)(b\ps)-\sigma(b\ps)(a\ps)=
(\sigma(a)\ps(b)-\sigma(b)\ps(a))\cdot\ps,
\end{multline} where $a,b\in\A$ and $\sigma$ is extended to a map
on $\A\cdot\ps$ by $\sigma(a\ps)=\sigma(a)\ps$. The elements
$e:=\ps, h:=-2t\ps$ and $f:=-t^2\ps$ span an $\F$-linear subspace
$$\mathcal{S}:=\linSpan_\F\{\ps,-2t\ps,-t^2\ps\}=\linSpan_\F\{e,h,f\}$$
of $\A\cdot\ps$. We restrict the multiplication (\ref{eq:bracket}) to
$\mathcal{S}$ without, at
this point, assuming closure. Now, $$\ps (t^2)=\ps(t\cdot
t)=\sigma(t)\ps(t)+\ps(t)t=(\sigma(t)+t)\ps(t)$$ which by using
(\ref{eq:bracket}), leads to
\begin{subequations}
\begin{align}
    \langle h,f\rangle&=2\langle
    t\ps,t^2\ps\rangle=2\sigma(t)\ps(t) t\ps,\label{eq:Shf}\\
    \langle h,e\rangle &= -2\langle t\ps,\ps\rangle =-2
    (\sigma(t)\ps(1)-\sigma(1)\ps(t))\ps\label{eq:She},\\
    \langle e,f\rangle &=-\langle
    \ps,t^2\ps\rangle=-(\sigma(1)(\sigma(t)+t)\ps(t)-
    \sigma(t)^2\ps(1)) \ps\label{eq:Sef}.
\end{align}
\end{subequations} Also, one would be tempted to make the general
ansatz
\begin{align*}
    \ps(1)=d_0+d_1t+\dots+d_kt^k,\qquad \sigma(1)=s_0+s_1t+\dots+s_lt^l.
\end{align*}If all non-negative integer powers of $t$ are linearly
independent over $\F$, then $\sigma(1)=1$ or $\sigma(1)=0$ since
$\sigma(1)=\sigma(1\cdot 1)=\sigma(1)^2$ and so $s_0$ is either $1$ or $0$
and $$s_1=\dots=s_l=0.$$ In addition to this, in this case,
\begin{multline*}
    d_0+d_1t+\dots+d_kt^k=\ps(1)=\ps(1\cdot
    1)=\sigma(1)\ps(1)+\ps(1)1=\\
    =(\sigma(1)+1)\ps(1)
    =(s_0+1)(d_0+d_1t+\dots+d_kt^k)
\end{multline*}leading to $d_0=d_1=\dots=d_k=0$ if $s_0=1$ and
arbitrary $d_0,\dots, d_k$ if $s_0=0$. But if $s_0=0$, that is
$\sigma(1)=0$, then $\sigma(t^w)=0$ for all $w\in\N$ since
$$\sigma(t^w)=\sigma(1\cdot t^w)=\sigma(1)\sigma(t^w)=0.$$ Under the
assumptions $\sigma(1)=1$ and $\ps(1)=0$ relations (\ref{eq:Shf}),
(\ref{eq:She}) and (\ref{eq:Sef}) simplify to
\begin{subequations}
\begin{align}
  &\langle h,f\rangle=2\sigma(t)\ps(t) t\ps\label{eq:Shfsimp}\\
  & \langle h,e\rangle =2\ps(t)\ps\label{eq:Shesimp}\\
  & \langle e,f\rangle =-(\sigma(t)+t)\ps(t)\ps\label{eq:Sefsimp}.
\end{align}
\end{subequations}
\begin{remark}
    Note that when $\sigma=\id$ and $\ps(t)=1$ we retain the
    classical $\sll$ with relations {\rm (\ref{eq:sl2})}.
\end{remark}
    \subsection{Quasi-Deformations with base algebra
    $\A=\F[t]$}\label{sec:k[t]}
Take $\A$ to be the polynomial algebra $\F[t]$, $\sigma(1)=1$ and
$\ps(1)=0$. Since all non-negative integer powers of $t$ are
linearly independent over $\F$ in $\F[t]$, we are in the situation
of relations (\ref{eq:Shf}), (\ref{eq:She}) and (\ref{eq:Sef}).
Suppose that $\sigma(t)=q(t)$ and $\ps(t)=p(t)$, where
$p(t),q(t)\in\F[t]$. To have closure of (\ref{eq:Shfsimp}),
(\ref{eq:Shesimp}) and (\ref{eq:Sefsimp}) these polynomials are
far from arbitrary. Indeed, by (\ref{eq:Shfsimp}) we get
$$\deg(\ps(t)\sigma(t))=\deg(p(t)q(t))\leq 1.$$ So three cases
arise
\begin{description}
    \item[\emph{Case 1}:]
    $\sigma(t)=q(t)=q_0+q_1t$, $q_1\neq 0$, $p(t)=p_0,$
    \item[\emph{Case 2}:] $\sigma(t)=q(t)=q_0$, $q_0\neq
    0$, $p(t)=p_0+p_1t,$
    \item[\emph{Case 3}:]
    $\sigma(t)=q(t)=0$, $p(t)=p_0+p_1t+\dots+p_nt^n$,
\end{description}
where in all three cases we assume $p(t)\neq 0$. Note that if
$p(t)=0$ then $\ps=0$ and so the original operator representation
collapses.
\begin{remark}
    If we allow $\sigma(t)=q(t)$ and $\ps(t)=p(t)$ where $p, q$
    are arbitrary polynomials in $t$ then we get a deformation of
    $\sll$ which does not preserve dimension; that is,
    brackets of the basis elements $e,f,h$ are not simply linear combinations in
    these elements but include more new "basis" elements. This phenomena could
    possibly be interesting to study further.
\end{remark}
\noindent\textbf{\textit{Case 1}:} Assume $q(t)=q_0+q_1t$,
implying that $p(t)=p_0$. Relations (\ref{eq:Shfsimp}),
(\ref{eq:Shesimp}) and (\ref{eq:Sefsimp}) according to
(\ref{eq:bracket}) now become
\begin{subequations}
    \begin{align}
        &\langle h,f\rangle :\,-2q_0ef+q_1hf+q_0^2eh-q_0q_1h^2-q_1^2fh=-q_0p_0h-2q_1p_0f
       \label{eq:Jackhf}\\
        &\langle h,e\rangle :\, -2q_0e^2+q_1he-eh= 2p_0e\label{eq:Jackhe}\\
        &\langle e,f\rangle :\, ef+q_0^2e^2-q_0q_1he-q_1^2fe=-q_0p_0e+\frac{q_1+1}{2}p_0h\label{eq:Jackef}.
    \end{align}
\end{subequations}
\begin{remark}Notice that changing the role of $h$ and $f$ in {\rm
(\ref{eq:Jackhf})} does not correspond to changing $h$ and $f$ in
{\rm (\ref{eq:bracket})}. This means that, in a sense, the
skew-symmetry of {\rm (\ref{eq:bracket})} is "hidden" in {\rm
(\ref{eq:Jackhf})}.  If $\langle f,h\rangle$ is calculated from
{\rm (\ref{eq:bracket})} one sees that indeed $\langle f,h\rangle
=-\langle h,f\rangle$ as one would expect and one gets exactly
minus the left-hand-side of {\rm (\ref{eq:Jackhf})}.
\end{remark}
Henceforth, we denote by $\mathbb{F}\{x_1,\dots, x_n\}$ the free
associative algebra over $\mathbb{F}$ on the set $\{x_1,\dots,
x_n\}$, i.e., the non-commutative polynomial (tensor) algebra over
$\mathbb{F}$ in the indeterminates $x_1,\dots, x_n$.

The associative algebra with three abstract generators $e,h$ and
$f$ and defining relations (\ref{eq:Jackhf}), (\ref{eq:Jackhe})
and (\ref{eq:Jackef}), that is, $\mathbb{F}\{e,f,h\}$ modulo the
relations defined by (\ref{eq:Jackhf}), (\ref{eq:Jackhe}) and
(\ref{eq:Jackef}), can be seen as a multi-parameter deformation of
$$\mathbf{U}(\sll)=\mathbb{F}\{e,f,h\}\big/
\begin{pmatrix}
[h,e]-2e,\quad [h,f]+2f,\quad [e,f]-h
\end{pmatrix},$$ where $[\cdot,\cdot]$ denotes the commutator.
\begin{example} By taking $q_0=0$ and $q:=q_1\neq 0$ we obtain a
deformation of $\sll$ corresponding to replacing the ordinary
derivation operator with the Jackson $q$-derivative $D_q$ (see
Section \ref{sec:qhlsigma}). Since this deformation therefore is
quite interesting we explicitly record it
\begin{align}\label{eq:Jacks_def}
        hf-qfh&=-2p_0f\notag\\
        he-q^{-1}eh&= 2q^{-1}p_0e\\
        ef-q^2fe&=\frac{q+1}{2}p_0h\notag.
\end{align}
 Note that, by taking $p_0=q=1$, we obtain the usual
commutation relations for $\sll$, corresponding to the ordinary
derivation operator $\partial$. We denote the "lifting" of the
right-hand-side of (\ref{eq:Jackhf}), (\ref{eq:Jackhe}),
(\ref{eq:Jackef}) for $p_0=1$, $q_0=0$, to an abstract
skew-symmetric algebra with products
$$\langle h,f\rangle = -2qf, \qquad \langle h,e\rangle =
2e,\qquad \langle e,f\rangle =\frac{q+1}{2}h,$$ by $\Jacksll$ and
call it
  informally the "Jackson
$\sll$". In Example
  \ref{exam:JackhomLie} we will see that $\Jacksll£$ is a (quasi-)
  hom-Lie algebra. When $p_0\neq 1$ we get a natural one-parameter
deformation of $\Jacksll$. The algebra $\F\{
  e,f,h\}/(\ref{eq:Jacks_def})$, with $p_0=1$ in
  (\ref{eq:Jacks_def}), can be thought of as an
 analogue of the universal enveloping algebra for $\Jacksll$. For
 $q=1$ it is indeed the universal enveloping algebra of the Lie
 algebra $\sll$. We denote by $\mathcal{U}_q$ the algebra  $\F\{
  e,f,h\}/(\ref{eq:Jacks_def})$. When $p_0\neq 1$ we similarly
  get a one-parameter deformation of $\mathcal{U}_q$.

There is also a Casimir-like element in this algebra, namely,
$$\Omega_q:=ef+qfe+\frac{q+1}{4}h^2=ef+qfe+\frac{\{2\}_q}{4}h^2.$$ For $q=1$ we retain the
classical central Casimir element for $\mathbf{U}(\sll)$ in the
basis $e,f,h$. It is straightforward to check that $\Omega_q$ is
normal in the sense that $\tau(z)\cdot\Omega_q=\Omega_q\cdot z$
for some map $\tau$ (notice that in the Lie case $q=1$ the element
$\Omega_1$ is really central, i.e., $\tau=\id$). The associative
product in $\mathcal{U}_q$ forces $\tau$ to be an algebra
endomorphism. Indeed, on the one hand: $\Omega_q\cdot
zw=\tau(z)\cdot\Omega_q\cdot w=\tau(z)\tau(w)\cdot\Omega_q$, and
on the other: $\Omega_q\cdot zw=\tau(zw)\cdot\Omega_q$. It can
also be checked that $\tau$ is in fact an automorphism. In our
present setting $\tau$ is completely determined by its action on
the basis elements as $\tau(e)=q^{-2}e$, $\tau(h)=h$ and
$\tau(f)=q^{2}f$. \qed
\end{example}
\begin{remark}
Assume that $q\neq 1$ and rewrite the relations {\rm
(\ref{eq:Jacks_def})} in the form{\rm :}
\begin{align}\label{eq:BSrelJacksl2}
\begin{split}
    yz-q^{-1}zy&=0\\
     zx-q^2xz&=p_0^2(1+q^{-1})y+p_0^2\frac{q+1}{q-1}\mathbf{1}\\
     xy-q^{-1}yx&=0,
\end{split}
\end{align}where we made the transformations $h\mapsto
2p_0q^{-1}y-2p_0/(q-1)\mathbf{1}, f\mapsto x$ and $e\mapsto z.$
Since the transformation is a bijection {\rm(}it is simply a
change of basis{\rm)} we deduce that
$\F\{e,f,h\}/(\ref{eq:Jacks_def})$ and $\F\{
x,y,z\}/(\ref{eq:BSrelJacksl2})$ are isomorphic as algebras. We
denote the algebra $\F\{ x,y,z\}/(\ref{eq:BSrelJacksl2})$ by
$\mathcal{W}_q$. Hence $\mathcal{W}_q\cong\mathcal{U}_q$. It is
easy to see that $\mathcal{W}_q$ is an iterated Ore extension as
follows. Indeed, putting
$$\mathcal{B}:=\mathbb{F}\{y,z\}\big/(yz-q^{-1}zy)$$ we see that
this is an Ore extension of $\F[z]$ with automorphism $z\mapsto
q^{-1}z$. Then extending once more we get $\mathcal{W}_q$ as the
Ore extension $\mathcal{B}[x,\varsigma,\partial_\varsigma]$ of
$\mathcal{B}$ defined by $\varsigma(z)=q^{-2}z$,
$\partial_\varsigma(z)=-q^{-2}ay-q^{-2}b\mathbf{1}$,
$\varsigma(y)=q^{-1}y$ and $\partial_\varsigma(y)=0$. {\rm(}What
is needed to check is that $\partial_\varsigma(yz)$ and
$\partial_\varsigma(q^{-1}zy)$ both give the same result.{\rm)} We
now note that iterated Ore extensions of an Auslander-regular
algebra {\rm(}see below for the definition{\rm)}, in this case
$\mathbb{F}[z]$, are themselves Auslander-regular by a theorem of
Ekstr\"om {\rm\cite{EkstromD}}. Moreover, one can prove
{\rm(\cite{BellSmithD}}, Proposition 2.1.1{\rm)} that
$\mathcal{W}_q$ has global dimension at most three. Also, using
Proposition 2.1.2 in {\rm\cite{BellSmithD}} or the Diamond Lemma
{\rm\cite{HellD}} it is easy to check that $\mathcal{W}_q$ has a
PBW-basis, hence is a noetherian domain of Gel'fand--Kirillov
dimension three {\rm(\cite{BellSmithD}}, Proposition 2.1.1{\rm)}.
Having a PBW-basis ensures Koszulity as an almost quadratic
algebra {\rm\cite{PiontSilvD}}. We note the following interesting
special case of relations {\rm(\ref{eq:BSrelJacksl2})}.
\begin{example}\label{exam:colJacksl2}
When $q=-1$ the endomorphism $\sigma$ becomes $\sigma(t)=-t$ and
the Jackson $q$-derivative is thus given by $f\mapsto
(f(t)-f(-t))/2t.$ The defining relations (\ref{eq:BSrelJacksl2})
for $\mathcal{W}_{-1}$ then become
\begin{align}\label{eq:colrelJacksl2}
    yz+zy=0,\qquad zx-xz=0,\qquad xy+yx=0.
\end{align}This is in fact a color-commutative Lie algebra graded by
$\Z_2^2$. Indeed, suppose we have a vector space $V$ with
decomposition
$$V=V_{(0,0)}\oplus
V_{(1,0)}\oplus V_{(0,1)}\oplus V_{(1,1)}=\F\{0\}\oplus \F y\oplus
\F\{0\}\oplus \big(\F x\oplus \F z\big).$$ Taking the colored
commutator $[A,B]_{\text{col}}:=AB-(-1)^{(\deg(A),\deg(B))}BA,$
where we have the bilinear form $(\cdot,\cdot)$ defined by
$(\deg(A),\deg(B)):=\alpha_1\beta_1+\alpha_2\beta_2$ for
$\deg(A)=(\alpha_1,\alpha_2)$ and $\deg(B)=(\beta_1,\beta_2)$
yields the algebra defined by relations (\ref{eq:colrelJacksl2}).
\qed
\end{example}
\end{remark}
Homogenizing the relations (\ref{eq:Jackhf}), (\ref{eq:Jackhe})
and (\ref{eq:Jackef}) with respect to a central degree-one element
$\zeta$ leads to algebras reminiscent of the central extensions of
three-dimensional Artin--Schelter regular algebras (in particular
the Sklyanin algebra) studied (among other algebras) in
\cite{LeBruynSmithvdBerghD}. The Artin--Schelter regular algebras
with homogeneous defining relations (and more generally
Auslander-regular algebras) turned out to be a most natural choice
for a non-commutative projective geometry and this stimulated a
directed effort in understanding the various homological and
geometric properties of these algebras (see \cite{ArtinSchelterD},
\cite{ArtinTatevdBerghD}, \cite{LeBruynSmithvdBerghD}, for
instance). Also, let $\g$ be a finite-dimensional Lie algebra.
Then the universal enveloping algebra $\mathbf{U}(\g)$ is filtered
by the canonical degree-filtration. This is a Zariski-filtration
and since $\mathrm{gr}(\mathbf{U}(\g))$ is isomorphic to a
commutative polynomial algebra (by the PBW-theorem) this being
Auslander-regular, it follows that $\mathbf{U}(\g)$ is itself
Auslander-regular \cite{LiOys}. Le Bruyn and Smith, Le Bruyn and
Van den Bergh, showed that $\mathbf{U}(\g)$ (or rather its
homogenizations) is connected to non-commutative projective
algebraic geometry \cite{LeBruynSmithD}, \cite{LeBruynvdBerghD}.

For the reader's convenience we recall the definition of
Artin--Schelter regular and Auslander-regular algebras. Let $R$ be
a graded connected algebra over a field $\F$, i.e.,
$R=\bigoplus_{n\in\Z_{\geq
    0}}R_n$ with $R_0\cong\F$. Then $R$
is Artin--Schelter regular (AS-regular) if it has finite global
dimension $d$ (all graded $R$-modules have finite projective
dimension $\leq d$) and
\begin{itemize}
\item $R$ has finite Gel'fand--Kirillov dimension, that is,
  if it has polynomial growth, i.e., if $\dim_\F(R_n)\leq n^k$, for some
  constant $k\in\N$;
\item $R$ is Gorenstein:
  $\mathrm{Ext}^i_{R}(\F,R)=0$ if $i\neq d$ and
  $\mathrm{Ext}^d_{R}(\F,R)\cong \F$.
\end{itemize}A related concept is the Auslander-regular
algebras. Let $R$ be a noetherian ring. Then $R$ is said to be
\emph{Auslander-regular} if for every $R$-module $\mathcal{M}$ and
for all
  submodules $\mathcal{N}$ of
  $\mathrm{Ext}^i_{R}(\mathcal{M},R)$ we have
   $\mathrm{min}\{j\mid
\mathrm{Ext}^j_{R}(\mathcal{N},R)\neq 0\}\geq i$, $\forall\,\,
i\geq 0$, in addition to $R$ having finite injective and global
dimension. For graded algebras of Gel'fand--Kirillov dimension
less than or equal to three it has been shown by Thierry Levasseur
\cite{LevasseurD} that Auslander-regular and AS-regular are in
fact equivalent.

We now review the notion of conformal $\sll$ enveloping algebras
introduced by Lieven Le Bruyn in \cite{LeBruynD}, showing that
$\mathcal{U}_q$ defined by (\ref{eq:Jacks_def}) also is a special
case of his construction.
\begin{example}\label{exam:LeBruynCEA}
Let $\mathbf{a}$ denote the $\C$-vector
  $(a_1,a_2,a_3,a_4,a_5,a_6,a_7)$. Define an algebra
$$\mathcal{F}(\mathbf{a}):=\C\{ x,y,z\}\Big /
\begin{pmatrix}
xy-a_1yx-a_2y,\\
yz-a_3zy-a_4x^2-a_5x,\\
zx-a_6xz-a_7z
\end{pmatrix}.$$ Notice that, for instance,
$\mathcal{F}(1,2,1,0,1,1,2)\cong\mathbf{U}(\mathfrak{sl}_2(\C))$.
Therefore it is tempting to view $\mathcal{F}(\mathbf{a})$, for
generic $\mathbf{a}$'s, as deformations of
$\mathbf{U}(\mathfrak{sl}_2(\C))$. Such a deformation of
$\mathbf{U}(\mathfrak{sl}_2(\C))$ is called a \emph{conformal
$\mathfrak{sl}_2$ enveloping algebra} if
$\mathrm{gr}(\mathcal{F}(\mathbf{a}))$, with
$\mathcal{F}(\mathbf{a})$ taken with the canonical
degree-filtration, is a three-dimensional Auslander-regular
quadr\-atic algebra. Edward Witten introduced a special case of
the above, namely the algebra
$$\mathsf{W}_a:=\mathcal{F}(a,1,1,a-1,1,a,1)=\C\{ x,y,z\}\Big /
\begin{pmatrix}
xy-ayx-y,\\
yz-zy-(a-1)x^2-x,\\
zx-axz-z
\end{pmatrix}$$related to some aspects of two-dimensional conformal field
theory. The case $a=1$ gives an algebra isomorphic to
$\mathbf{U}(\mathfrak{sl}_{2}(\C))$.

Suppose $a_1a_2a_3a_5a_6a_7\neq 0$. Le Bruyn shows that under this
condition $\mathcal{F}(\mathbf{a})$ is a conformal
$\mathfrak{sl}_2(\C)$ enveloping algebra if and only if $a_6=a_1$
and $a_7=a_2$.

We now note that when $\F=\C$ our algebra $\mathcal{U}_q$ is a
conformal $\mathfrak{sl}_2(\C)$ enveloping algebra. In fact
$\mathcal{U}_q=\mathcal{F}(\mathbf{a'})$, under the assignments
$h\mapsto x$, $e\mapsto y$, $f\mapsto z$ and where
$$\mathbf{a'}=(q^{-1},2q^{-1}p_0,q^2,0,\frac{q+1}{2}p_0,q^{-1},2q^{-1}p_0).
$$This also shows that $\mathcal{U}_q$ is Auslander-regular since $a_1=a_6$ and
$a_2=a_7$. \mbox{Le Bruyn} also shows in \cite{LeBruynD} that
there is a very nice non-commutative geometry behind these
algebras when homogenized.\qed
\end{example}
By taking $p_0=q_0=0$ in (\ref{eq:Jacks_def}) we get an
``abelianized'' version of $\Jacksll$. However, as we mentioned
before, the operator representation we started with collapses in
this case since $\ps=0$.
\begin{remark} Note that, unlike the undeformed $\mathbf{U}(\sll)$, the
  one-parameter deformation of $\,\mathcal{U}_q$
  {\rm(}with $q\neq 1$, $p_0\neq 0${\rm)} has non-trivial one-dimensional representations. Indeed, taking
  $h={2p_0}/({q-1})$ we see from the defining relations that $ef=(q+1)p_0^2/(q-1)$ and so $e$ and
  $f$ can be chosen as any numbers satisfying this equality. When
  $q\to 1$ in which case we ``approach $\sll$'' we see that
  $h$ goes to infinity and so the representation collapses in the limit.
\end{remark}
\begin{remark}
It is important to notice that Theorem \ref{thm:GenWitt} gives us two
alternative
ways to view our algebras{\rm :}
\begin{itemize}
     \item[{\rm (i)}] either one uses Theorem \ref{thm:GenWitt} to
       construct a bracket on
       the vector space $\A\cdot\ps$ with the help of equation
       {\rm (\ref{eq:GenWittProdFormula})}, viewing this bracket as a product;
     \item[{\rm (ii)}]  or one uses {\rm
         (\ref{eq:GenWittProdFormula})} in adjunction
       with the definition {\rm (\ref{eq:GenWittProdDef})} of the bracket to
       obtain quadratic relations and ``moding''
       out these from the tensor algebra, thereby giving us
       an associative quadratic algebra which loosely
       can be thought of as a deformed ``enveloping algebra'' of the
       corresponding algebra from viewpoint {\rm (i)}, or even as a deformed
       universal enveloping algebra for $\sll$.
\end{itemize}
\end{remark}
\subsubsection*{Twisted Jacobi identity for Case 1}
The possible $\delta\in\A=\F[t]$ in the twisted Jacobi identity
can be computed from the condition
\begin{align}\label{eq:deltacond}
\ps\circ\sigma(a)=\delta\cdot\sigma\circ\ps(a),
\end{align}
required to hold for any $a\in\A$. For $a=t^0=1$ this equation
holds trivially because $\sigma(1)=1$ and $\ps(1)=0$. For $a=t$
the left-hand-side becomes
\begin{align*}
    \ps\circ\sigma(t)=\ps(q_0+q_1t)=q_1p_0
\end{align*}and the right-hand-side becomes
\begin{align*}
    \delta\cdot\sigma\circ\ps(t)=\delta\cdot\sigma(p_0)=\delta\cdot
    p_0,
\end{align*}from which we immediately see that $\delta=q_1$. To show
that (\ref{eq:deltacond}) holds for all $a\in\A$ with $\delta=q_1$
it is enough by linearity to show this for arbitrary monomial
$a=t^k$. As we have seen the statement is true for $k=0,1$. Assume
the statement holds for $k=l$. On writing $t^{l+1}=tt^l$ and using
the $\sigma$-Leibniz rule (\ref{eq:sigmaLeibniz}) we get:
\begin{multline*}
     \mathrm{R.H.S}=\delta\cdot\sigma\circ\ps(t^{l+1})=\delta\cdot
     (\sigma\circ\ps(t)\sigma(t)^l+\sigma^2(t)\sigma\circ\ps(t^l)
     )=\\
     =\delta\cdot\sigma\circ\ps(t)\sigma(t^l)+\sigma^2(t)\cdot\delta\cdot\sigma\circ\ps(t^l)=[\text{induction
     step}]=\\
     =(\ps\circ\sigma(t)
     )\sigma(t^l)+\sigma^2(t)\ps\circ\sigma(t^l)=
     \ps\circ\sigma(t^{l+1})=\mathrm{L.H.S}.
\end{multline*}Note that this shows that it is enough to check
condition (\ref{eq:deltacond}) for low degrees, even in cases
other than the present. So by Theorem \ref{thm:GenWitt} we now have a
deformed Jacobi identity
$$\cs_{x,y,z}(\langle \sigma(x),\langle
y,z\rangle\rangle+q_1\langle x,\langle y,z\rangle\rangle)=0$$ on
$\A\cdot\ps=\F[t]\cdot\ps$. By defining
$\alpha(x):=q_1^{-1}\sigma(x)$ this can be re-written as the
Jacobi-like relation for a hom-Lie algebra
\begin{equation}\label{eq:JacCase1}
\cs_{x,y,z}\langle (\alpha+\id)(x),\langle y,z\rangle\rangle=0.
\end{equation}
Note that this twisted Jacobi identity follows directly from the
general Theorem \ref{thm:GenWitt}. This shows that we do not have
to make any lengthy trial computations with different more or less
ad hoc assumptions on the form of the identity.
\begin{example}\label{exam:JackhomLie}
The algebra $\Jacksll$ is a hom-Lie algebra with
$\alpha=q^{-1}\sigma$ and with Jacobi identity as given by
(\ref{eq:JacCase1}). This is not clear \`a priori since $\Jacksll$
is \emph{not} a subalgebra of $\mathcal{A}\cdot\ps$. When
$\Jacksll$ is represented by its "non-lifted" operators $e=\ps,
h=-2t\ps$ and $f=-t^2\ps$, this representation becomes a
subalgebra of $\mathcal{A}\cdot\ps$ and hence a hom-Lie algebra.
However, $\Jacksll$ is defined abstractly without reference to a
particular representation. It is initially possible that the
twisted Jacobi identity is a result of some extra relations (e.g.,
by the $\sigma$-Leibniz rule) available in that representation as
$\sigma$-derivations.

To prove that $\Jacksll$ is a hom-Lie algebra we proceed as
follows. First define $\alpha$ on the basis vectors $e,f,h$ as
$\alpha(e)=q^{-1}e$, $\alpha(h)=h$ and $\alpha(f)=qf$. It is clear
that it is enough to consider the case when $x=e$, $y=f$ and
$z=h$. Then (\ref{eq:JacCase1}) becomes
\begin{align*}
\langle (\alpha+\id)(e),&\langle f,h\rangle\rangle+\langle
(\alpha+\id)(f),\langle h,e\rangle\rangle+\langle
(\alpha+\id)(h),\langle e,f\rangle\rangle=\\
&=2q(q^{-1}+1)\langle e,f\rangle+2(q+1)\langle
f,e\rangle+(q+1)\langle h,h\rangle =\\
&=(2(q+1)-2(q+1))\langle e,f\rangle=0,
\end{align*}where we have used that $\langle\cdot,\cdot\rangle$ is
skew-symmetric.\qed\end{example}
 \noindent\textbf{\textit{Case 2}:} Assume now that $\sigma(t)=q(t)=q_0\neq 0$
and so $\ps(t)=p(t)=p_0+p_1t$. The relations (\ref{eq:Shfsimp}),
(\ref{eq:Shesimp}) and (\ref{eq:Sefsimp}) combined with
(\ref{eq:bracket}) become
\begin{subequations}
    \begin{align*}
       &\langle h,f\rangle :\quad -2ef+q_0eh=-p_0h-2p_1f\\
       &\langle h,e\rangle :\quad -2q_0e^2-eh= 2p_0e-p_1h\\
       &\langle e,f\rangle :\quad
       ef+q_0^2e^2=-q_0p_0e+\frac{q_0p_1+p_0}{2}h+p_1f,
    \end{align*}
\end{subequations}since $\sigma(e)=e$, $\sigma(h)=-2q_0e$ and
$\sigma(f)=-q_0^2e.$
    It is obvious that we cannot recover the classical $\sll$ from this
    deformation by specifying suitable parameters since, in a sense,
    the commutators "collapsed" due to the choice of $\sigma$.
\subsubsection*{Twisted Jacobi identity for Case 2}
It is easy to deduce, following the same arguments as in
\emph{Case 1}, that one can take $\delta=0$ in this case and so
have \begin{align}\label{eq:JacCase2}
\cs_{x,y,z}\langle
\sigma(x),\langle y,z\rangle
\rangle=0
\end{align} on $\A\cdot\ps$ making it into a hom-Lie algebra.\\

\noindent\textbf{\textit{Case 3}:} It follows immediately upon
 insertion of $\ps(t)=p(t)=p_0+p_1t+\dots+p_nt^n$
 into (\ref{eq:Sefsimp}) that $\deg(p(t))\leq 1$. With this
 (\ref{eq:Shfsimp}),
(\ref{eq:Shesimp}) and (\ref{eq:Sefsimp}) combined with
(\ref{eq:bracket}) now become
\begin{align*}
  0=0,\qquad -eh=
    2p_0e-p_1h, \qquad
       ef=\frac{p_0}{2}h+p_1f.
\end{align*}
\subsubsection*{Twisted Jacobi identity for Case 3}
In this case (\ref{eq:deltacond}) becomes $0=\ps(0)=\ps(
\sigma(t^k))=\delta\cdot \sigma(\ps(t^k))=\delta\cdot 0$, valid
for all $k\in\Z_{\geq 0}$, and so $\delta$ can be chosen
completely arbitrary in this case, for example $\delta=0$ giving
again (\ref{eq:JacCase2}) and thus a hom-Lie algebra.
\subsection{Deformations with base algebra $\F[t]/(t^3)$}\label{sec:k[t]/t^3}
Now, let $\F$ include all third roots of unity and take as $\A$
the algebra $\F[t]/(t^3)$. This is obviously a three-dimensional
$\F$-vector space and a finitely generated $\F[t]$-module with
basis $\{1,t,t^2\}$. Let, as before, $e=\ps, h=-2t\ps,$ and
$f=-t^2\ps.$ Note that $-2t\cdot e=h$, $t\cdot h=2f$ and $t\cdot
f=0$. Put
\begin{subequations}
\begin{align}
    \ps(t)&=p_0+p_1t+p_2t^2\label{eq:psz}\\
    \sigma(t)&=q_0+q_1t+q_2t^2\label{eq:sigmaz}
\end{align}
\end{subequations}considering these as elements in the ring $\F[t]/(t^3)$.
We will once again make the assumptions that $\sigma(1)=1$,
$\ps(1)=0$ and so relations (\ref{eq:Shfsimp}), (\ref{eq:Shesimp})
and (\ref{eq:Sefsimp}) still hold. The equalities (\ref{eq:psz})
and (\ref{eq:sigmaz}) have to be compatible with $t^3=0$. This
means in particular that
\begin{align*}
    0=\sigma(t)^3=(q_0+q_1t+q_2t^2)^3=q_0^3+3q_0^2q_1t+(3q_0^2q_2+3q_0q_1^2)t^2
\end{align*}implying $q_0=0$. Moreover,
\begin{align}\label{eq:psz3}
0=\ps(t^3)=\sigma(t)^2\ps(t)+\ps(t^2)t=
(\sigma(t)^2+(\sigma(t)+t)t)\ps(t).
\end{align}
The action of $\A\cdot\ps$ on $\A$ has a matrix
representation in the basis $\{1,t,t^2\}$:
\begin{align*}
    &e\mapsto
    \begin{pmatrix}
        0& p_0 & 0\\
        0 &p_1&(q_1+1)p_0\\
        0 &p_2&(q_1+1)p_1+q_2p_0
    \end{pmatrix},
    &h\mapsto
    \begin{pmatrix}
        0& 0 & 0\\
        0 &-2p_0&0\\
        0 &-2p_1&-2(q_1+1)p_0
    \end{pmatrix}
\end{align*}and
\begin{align*}
    &f\mapsto
    \begin{pmatrix}
        0& 0 & 0\\
        0 &0&0\\
        0 &-p_0&0
    \end{pmatrix}.
\end{align*}Note, however, that in the case when $p_0=1$,
$p_1=p_2=0$, $q_0=q_2=0$ and $q_1=1$, which seemingly would
correspond to the case of the classical $\sll$, we surprisingly
get the matrices
\begin{align*}
    &e\mapsto
    \begin{pmatrix}
        0& 1 & 0\\
        0 &0&2\\
        0 &0&0
    \end{pmatrix},
    &h\mapsto
    \begin{pmatrix}
        0& 0 & 0\\
        0 &-2&0\\
        0 &0&-4
    \end{pmatrix},\qquad
    &f\mapsto
    \begin{pmatrix}
        0& 0 & 0\\
        0 &0&0\\
        0 &-1&0
    \end{pmatrix}
\end{align*}satisfying the "commutation" relations
\begin{align*}
    hf-fh=-2f,\qquad
    he-eh=2e,\qquad
    ef+2fe=h.
\end{align*}The reason for this anomaly is that the values of the
parameters so chosen are not compatible with the $\sigma$-Leibniz
rule on $\A=\F[t]/(t^3)$ since (\ref{eq:psz3}) becomes $0=3t^2$
this equality not possible in $\A$ over a field of characteristic
zero.

The bracket can be computed abstractly on generators as
\begin{subequations}
\begin{align}
    \label{eq:brack_hf}&\langle
    h,f\rangle=q_1hf+2q_2f^2-q_1^2fh\\
    \label{eq:brack_he}&\langle h,e\rangle
    =q_1he+2q_2fe-eh\\
    \label{eq:brack_ef}&\langle
    e,f\rangle=ef-q^2_1fe.
\end{align}
\end{subequations}
 Formulas (\ref{eq:psz}) and
(\ref{eq:sigmaz}) together with the assumption that the
right-hand-side of these are elements in $\F[t]/(t^3)$ now yield
closure. Indeed, by (\ref{eq:bracket}), we get
\begin{subequations}
\begin{eqnarray}
    \langle h,f\rangle &=&
    2(q_1t+q_2t^2)(p_0+p_1t+p_2t^2)t\ps=\notag\\
    &=&2q_1p_0t^2\ps=-2q_1p_0f,\label{eq:hf}\\
    \langle h,e\rangle &=&2
    (p_0+p_1t
    +p_2t^2
    )\ps
    =2p_0e-p_1h
    -2p_2f,\label{eq:he}\\
    \langle e,f\rangle &=&((q_1+1)t+q_2t^2)(p_0+p_1t+p_2t^2)\ps=\notag\\
    &=& (p_{{1}}+q_{{1}}p_{{1}}+q_{{2}}p_{{0}}) f
 +\frac{q_1+1}{2}p_0 h. \label{eq:ef}
\end{eqnarray}
\end{subequations}
We have one other condition to take into
account, namely the relation (\ref{eq:psz3}). This relation gives
us that
$$(\sigma(t)^2+(\sigma(t)+t)t)\ps(t)=p_0(q_1^2+q_1+1)t^2=0.$$ In other
words, if $p_0\neq 0$, we
generate deformations at the imaginary third roots of unity; if
$p_0=0$ then $q_1$ is a true formal deformation parameter. This
means that we once again have to subdivide our presentation from
this point
into two cases.\\

\noindent \textbf{\textit{Case 1}:} When $p_0=0$ the relations for
the deformed $\sll$ become
\begin{align}\label{eq:rootunitycase1}
\begin{split}
    &\langle h,f\rangle :\quad q_1hf+2q_2f^2-q_1^2fh=0\\
    &\langle h,e\rangle :\quad q_1he+2q_2fe-eh=-p_1h
    -2p_2f\\
    &\langle e,f\rangle :\quad ef-q^2_1fe= p_1(q_{{1}}+1)
    f.
\end{split}
\end{align} Once again, we cannot recover $\sll$ from this
deformation by choosing suitable values of the parameters since
$\langle e,f\rangle$ can never be $h$, $\langle h,e\rangle$ can
never be $2e$ and $\langle h,f\rangle$ is identically zero.

Take $p_1=p_2=0$ (so $\ps$ is identically zero), $q_1=1$ and
$\epsilon:=-2q_2$ in (\ref{eq:rootunitycase1}). Then the first
relation in (\ref{eq:rootunitycase1}) becomes
\begin{align}\label{eq:JordQuant}
hf-fh=\epsilon f^2
\end{align} which is the defining
relation for the \emph{Jordanian quantum plane} or the
\emph{$\epsilon$-deformed quantum plane}. In \cite{ChoMadoreParkD} is
developed the differential geometry of the Jordanian quantum plane in
analogy with commutative geometry's moving frame formalism. This means
that they construct a suitable differential calculus (differential
1-forms). They also show that one can extend the Jordanian quantum plane
 by adjoining the inverses
to $h$ and $f$ and furthermore construct a differential calculus on
this space. It turns out that
the commutative limit of this extended Jordanian quantum plane has
a metric of constant Gaussian curvature $-1$ and so the authors argue that the
extended Jordanian quantum plane with this calculus is a deformation
of the Poincar\'e upper-half plane.
 The Jordanian and ordinary quantum
planes are both associated to (different) quantizations of the Lie
group $\mathrm{GL}_2(\F)$ in the sense that these quantizations are
symmetry groups with central determinants for the respective
plane. Considerations like these lead unavoidably into the realm of
quantum groups.

For $\epsilon=1$, the quadratic algebra with defining relation
(\ref{eq:JordQuant}) is actually
one of two possible Artin--Schelter regular algebras (graded in
degree one) of global dimension two, the other being the ordinary
quantum plane
$hf-qfh=0$ (see \cite{ArtinSchelterD}).

Suppose for this paragraph that $\F=\C$. If $q_1=1$, $p_1=1$,
$p_2=a/2$ and $q_2=0$ then (\ref{eq:rootunitycase1}) become the
one-parameter family of three-dimensional Lie algebras $\g$ with
relations:
\begin{align*}
    hf-fh=0,\qquad he-eh=-h-af,\qquad
    ef-fe=2f,\qquad\text{for } a\in\C.
\end{align*}
These algebras are solvable\footnote{In fact, it is easy to see
that $\g^{(2)}=0$ where we put $\g^{(1)}=\langle\g,\g\rangle$ and
then inductively defining $\g^{(i)}=\langle
\g^{(i-1)},\g^{(i-1)}\rangle$.} for every $a\in\C$. By the general
classification of three-dimensional complex Lie algebras this
means that $\g$ is isomorphic to one \cite{GorbOnishVinD} of the
Lie algebras (zero brackets are omitted)
\begin{itemize}
\item $\ell(2)\oplus\C e_3$: $\langle e_1,e_2\rangle=e_2$;
\item $\ell(3)$: $\langle e_1,e_2\rangle =e_2$, $\langle
  e_1,e_3\rangle =e_2+e_3$;
\item $\ell(3,c)$: $\langle e_1,e_2\rangle =e_2$, $\langle
  e_1,e_3\rangle =c\cdot
  e_3$, for $c\in\C$.
\end{itemize}
Furthermore, putting $p_1=0$
and $p_2=-1/2$ we get
\begin{align*}
    hf-fh=0,\qquad he-eh=f,\qquad ef-fe= 0,
\end{align*}the relations for the Heisenberg Lie algebra. Putting
instead $p_1=p_2=0$ we get the three-dimensional polynomial
algebra in three commuting variables. So, in a sense, the class of
algebras with three generators and defining relations
(\ref{eq:rootunitycase1}), being obtained via a special twisting
of $\sll$, is a multi-parameter deformation of the polynomial
algebra in three commuting variables, of the Heisenberg Lie
algebra and of the above solvable Lie algebra.
\subsubsection*{Twisted Jacobi identity for Case 1}
The left-hand-side of (\ref{eq:deltacond}) for $a=t$ is
\begin{multline*}
    \ps\sigma(t)=\ps(q_1t+q_2t^2)=q_1(p_1t+p_2t^2)+q_2(\sigma(t)+t)(p_1t+p_2t^2)=\\
    =q_1p_1t+q_1p_2t^2+q_2((q_1+1)t+q_2t^2)(p_1t+p_2t^2)=\\
=q_1p_1t+(q_1p_2+q_2p_1(q_1+1))t^2
\end{multline*} keeping in mind $t^3=0$. The right-hand-side becomes
\begin{multline*}
    \delta\cdot\sigma\circ\ps(t)=\delta\cdot(p_1t(q_1+q_2t)+p_2t^2(q_1+q_2t)^2)
    =\delta(p_1q_1t+(p_1q_2+p_2q_1^2)t^2).
\end{multline*}Assuming that $\delta$ can be written as $\delta=\delta_0+\delta_1
t+\delta_2t^2$ we let $\xi_0,\xi_1$ and $\xi_2$ denote three free
parameters. Then $\ps\circ\sigma(t)=\delta\cdot\sigma\circ\ps(t)$ becomes
$$q_1p_1t+(q_1p_2+q_2p_1(q_1+1))t^2=\delta_0p_1q_1t+(\delta_1p_1q_1+\delta_0(p_1q_2+p_2q_1^2))t^2$$
which is equivalent to the linear system of equations for
$\delta_0, \delta_1$ and $\delta_2$
\begin{align*}
\left\{
    \begin{array}{l}
        p_1q_1\delta_0-q_1p_1=0\\
        (p_1q_2+p_2q_1^2)\delta_0+p_1q_1\delta_1-q_2p_1(q_1+1)-q_1p_2=0
    \end{array}\right.
\end{align*}Since $\delta_2$ is not involved at all it can be
chosen arbitrary, say $\delta_2=\xi_2$.
 We then have several cases to consider:
\begin{enumerate}
    \item In the case $q_1p_1\neq 0$ we find
        $$(\delta_0, \delta_1,\delta_2)^{\mathrm{T}}=
            \big(1,
             -\frac{q_1p_2-p_2-p_1q_2}{p_1},
             \xi_2\big)^{\mathrm{T}}. $$ The twisted Jacobi identity
             thus becomes
        $$\cs_{x,y,z}\big (\langle \sigma(x),\langle
        y,z\rangle
        \rangle+(1-\frac{q_1p_2-p_2-p_1q_2}{p_1}t+\xi_2t^2)\langle
        x,\langle y,z\rangle\rangle\big)=0.$$
        Notice that this defines a quasi-hom-Lie algebra
        which is not a hom-Lie algebra.
    \item If $p_1\neq0, q_2\neq 0, q_1=0$  we get
        $\delta=1+\xi_1t+\xi_2t^2$.
        The twisted Jacobi identity can now be written as
        $$\cs_{x,y,z}\big (\langle \sigma(x),\langle
        y,z\rangle
        \rangle+(1+\xi_1t+\xi_2t^2)\langle
        x,\langle y,z\rangle\rangle\big )=0.$$ Since $\xi_1$
        and $\xi_2$ are arbitrary this equality is equivalent to
        \begin{align}\label{eq:JacCase12unity}
        \cs_{x,y,z}\langle (\sigma+\id)(x),\langle
        y,z\rangle\rangle=0
        \end{align} which means that we have a
        hom-Lie algebra for these parameters.
    \item $p_1\neq0$, $q_1=q_2=0$ yield
    $\delta=\xi_0+\xi_1t+\xi_2t^2$.
        The deformed Jacobi identity is
        $$\cs_{x,y,z}\big (\langle \sigma(x),\langle
        y,z\rangle
        \rangle+(\xi_0+\xi_1t+\xi_2t^2)\langle
        x,\langle y,z\rangle\rangle\big )=0$$
        and once again, since $\xi_1$
        and $\xi_2$ are arbitrary we get a hom-Lie algebra with the same twisted Jacobi
        identity as above (\ref{eq:JacCase12unity}).
    \item For $p_1=0,p_2\neq0,q_1\neq0$ we get
    $\delta=q_1^{-1}+\xi_1t+\xi_2t^2$.
            The twisted Jacobi identity is
        $$\cs_{x,y,z}\big (\langle \sigma(x),\langle
        y,z\rangle
        \rangle+(q_1^{-1}+\xi_1t+\xi_2t^2)\langle
       x,\langle y,z\rangle\rangle\big )=0.$$
        Since $\xi_1$ and $\xi_2$ are arbitrary they can be chosen to be zero and
        re-scaling $\sigma$ we get the deformed Jacobi identity of a hom-Lie algebra
        in this case as well. Lastly, we have
        \item $p_1=p_2=0$ or $p_1=q_1=0$ which both yield the same
        result as 3.
   \end{enumerate}
\noindent\textbf{\textit{Case 2}:} When $q_1=\omega^k\neq 1$,
where $\omega:=\mathrm{e}^{\frac{2\pi}{3}i}$ is a third root of unity and
$p_0\neq 0$ we have that the relations (\ref{eq:hf}),
(\ref{eq:he}) and (\ref{eq:ef}) become
\begin{align*}
    &\langle h,f\rangle :\quad
    \omega^khf+2q_2f^2-\omega^{2k}fh=-2\omega^kp_0f\\
    &\langle h,e\rangle :\quad \omega^khe+2q_2fe-eh=2p_0e-p_1h
    -2p_2f\\
    &\langle e,f\rangle :\quad ef-\omega^{2k}fe=((\omega^k+1)p_{{1}}+q_{{2}}p_{{0}}) f
 +\frac{\omega^k+1}{2}p_0 h,
\end{align*}where $k=1,2$.
Note that this is a deformation of the same type as
(\ref{eq:Jacks_def}) if specifying $p_1=p_2=q_2=0$.
\subsubsection*{Twisted Jacobi identity for Case 2}
The left-hand-side of (\ref{eq:deltacond}) for $a=t$ is
\begin{multline*}
    \ps\circ\sigma(t)=\ps(\omega^k t+q_2t^2)=\\
    =\omega^k(p_0+p_1t+p_2t^2)+
    q_2((\omega^k+1)t+q_2t^2)(p_0+p_1t+p_2t^2)=\\
    =\omega^k p_0+((p_1+p_0q_2)\omega^k+p_0q_2)t+((p_2+p_1q_2)\omega^k+(p_0+p_1)q_2)t^2
\end{multline*}and the right-hand-side is
\begin{multline*}
\delta\cdot\sigma\circ\ps(t)=\delta\cdot\sigma(p_0+p_1t+p_2t^2)
=\delta\cdot(p_0+p_1\omega^k t+(p_1q_2+p_2\omega^{2k})t^2).
\end{multline*} Assuming once again
$\delta=\delta_0+\delta_1t+\delta_2t^2$ and remembering
$\omega^3=1$, we get
$$\begin{pmatrix}
    \delta_0\\
    \delta_1\\
    \delta_2
  \end{pmatrix}=
  \begin{pmatrix}
   \omega^k\\
    \frac{-p_1\omega^{2k}+(p_1+q_2p_0)\omega^k+q_2p_0}{p_0}\\
    \frac{-(p_2p_0+p_1^2+p_1q_2p_0)\omega^{2k}+p_0(p_2-p_1q_2)\omega^k+p_1^2+q_2^2p_0^2+p_1q_2p_0}{p_0}
  \end{pmatrix}:=\begin{pmatrix}
    \omega^k\\
    \mathbf{w}_1\\
    \mathbf{w}_2
   \end{pmatrix}.
$$The deformed Jacobi identity can thus be written in the
following form
    $$\cs_{x,y,z}\big (\langle \sigma(x),\langle
        y,z\rangle
        \rangle
       +(\omega^k+\mathbf{w}_1t+\mathbf{w}_2t^2)\langle x,\langle
        y,z\rangle\rangle\big)=0,$$ this being the defining relation for a quasi-hom-Lie algebra.
This quasi-hom-Lie algebra becomes a hom-Lie algebra when
$\mathbf{w}_1=\mathbf{w}_2=0$, that is for some special choices of
parameters defining $\sigma(t)$ and $\ps(t)$.
\begin{remark}[Generating deformations at
$N^{\mathrm{th}}$-roots of unity] The construction above to
generate deformations at third roots of unity can be generalized.
Suppose $\A=\F[t]/(t^N)$, $N\geq 3$, and $p_0\neq 0$ and that $\F$
is a field which includes all {\rm(}primitive{\rm)}
$N^{\text{th}}$-roots of unity. Assume also that
$$\sigma(t)=t(q_1+q_2t)\qquad \text{and}\qquad
\ps(t)=p_0+p_1t+\dots+ p_{N-1}t^{N-1}.$$ The condition that
$q_0=0$ remains unaltered. We now have
\begin{multline*}
    0=\ps(t^N)=\sum_{j=0}^{N-1}\sigma(t)^jt^{N-j-1}\ps(t)=
    \sum_{j=0}^{N-1}t^j(q_1+q_2t)^jt^{N-j-1}\ps(t)=\\
    =t^{N-1}\sum_{j=0}^{N-1}(q_1+q_2t)^j(p_0+p_1t+\dots
    +p_{N-1}t^{N-1})
    =p_0\{N\}_{q_1}t^{N-1}
\end{multline*}keeping in mind $t^N=0$. This gives that $q_1$ is a $N^{\mathrm{th}}$-root
of unity. Once again taking $q_2=0$ and $p_1=p_2=\dots=p_{N-1}=0$
yields a deformation of type {\rm (\ref{eq:Jacks_def})} at
$N^{\mathrm{th}}$-roots of unity.
\end{remark}
\subsubsection*{Acknowledgments.}
We thank Petr Kulish, Gunnar Traustason, Gunnar Sigurdsson and
Freddy Van Oystaeyen for valuable comments.

This work has been partially supported by the Crafoord Foundation,
the Swedish Royal Academy of Sciences, the Mittag-Leffler
Institute and the Liegrits network.


\begin{thebibliography}{99}
\bibitem{ArtinSchelterD} Artin, M., Schelter, W., \emph{Graded Algebras
    of Global
    Dimension 3}, Adv. Math., \textbf{66} (1987), 171--216.
\bibitem{ArtinTatevdBerghD} Artin, M., Tate, J., Van den Bergh, M.,
  \emph{Some Algebras Associated to Automorphisms of Elliptic
  Curves},
  The Grothendieck Festschrift Vol 1, 33--85, Birkhauser, Boston (1990).
\bibitem{BellSmithD} Bell, A.D., Smith, S.P., \emph{Some 3-dimensional
    Skew Polynomial Rings}, Preprint, Draft of April 1, 1997.
\bibitem{ChaiKuLukD} Chaichian, M., Kulish, P.,
Lukierski, J., \emph{$q$-deformed Jacobi identity, $q$-oscillators
and $q$-deformed infinite-dimensional algebras}, Phys. Lett. B \textbf{237}
(1990), no. 3-4, 401--406.
\bibitem{ChoMadoreParkD} Cho, S., Madore, J., Park, K.S.,
  \emph{Noncommutative Geometry of the $h$-deformed Quantum Plane},
  Preprint, (1997), \texttt{q-alg/9709007}.
\bibitem{CurtrZachos1D} Curtright, T. L., Zachos, C. K.,
\emph{Deforming maps for quantum algebras}, Phys. Lett. B \textbf{243}
(1990), no. 3, 237--244.--478.
\bibitem{DamKuD} Damaskinsky, E. V., Kulish, P. P., \emph{Deformed
oscillators and their applications}
(Russian),  Zap.\ Nauch.\ Semin.\ LOMI
\textbf{189} 1991, 37--74; Engl. transl. in  J. Soviet Math.  \textbf{62}
no.~5, (1992), 2963--2986.
\bibitem{QuantFieldStringsD} Deligne, P., Etingof, P., Freed, D.S.,
Jeffrey, L.C., Kazhdan, D., \mbox{Morgan, J.W.,} Morrison, D.R.,
Witten, E., (Eds) \emph{Quantum Fields and Strings: A Course for
Mathematicians}, 2 vol., Amer. Math. Soc., 1999.
\bibitem{DiFranMathiSenCD} Di
Francesco, P., Mathieu, P., S\'en\'echal, D., \emph{Conformal
Field Theory}, Springer Verlag, 1997, 890 pp.
\bibitem{Drin1D} Drinfel'd, V.G., \emph{Hopf algebras and the quantum
Yang--Baxter equation}, Soviet Math. Doklady \textbf{32} (1985), 254--258.
\bibitem{EkstromD} Ekstr\"om, E.K., \emph{The Auslander condition
on graded and filtered noetherian rings}, S\'eminaire
Dubreil--Malliavin 1987-88, LNM 1404, Springer-Verlag, 1989,
220--245.
\bibitem{FLMD} Frenkel, I., Lepowsky, J., Meurman, A., \emph{Vertex
    Operator Algebras and the Monster},
    Academic Press, 1988, 508 pp.
\bibitem{Fuchs1D} Fuchs, J., \emph{Affine Lie Algebras and Quantum
Groups}, Cambridge University Press, 1992, 433 pp.
\bibitem{Fuchs2D} Fuchs, J., \emph{Lectures on Conformal Field
Theory and Kac-Moody Algebras}, Springer Lecture Notes in Physics 498,
(1997), 1--54.
\bibitem{GorbOnishVinD} Gorbatsevich, V.V., Onishchik, A.L.,
  Vinberg, E.B., \emph{Structure of Lie Groups and Lie Algebras} in
  Encyclopedia of Math. Sciences Vol. 41, Springer-Verlag 1994, 248 pp.
\bibitem{HartLarsSilv1D} Hartwig, J.T., Larsson, D., Silvestrov, S.D.,
  \emph{Deformations of Lie algebras using $\sigma$-derivations}, Preprints
  in Mathematical Sciences 2003:32, LUTFMA-5036-2003, Centre for
  Mathematical Sciences, Department of Mathematics, Lund Institute of
  Technology, Lund University, (2003), \texttt{math.QA/0408064},
  to appear in Journal of Algebra.
\bibitem{HellD} Hellstr\"om, L., \emph{The Diamond Lemma for Power
Series Algebras}, Doctoral Thesis, no 23, 2002, Ume\aa\, University.
\bibitem{HelSil-bookD}{\rm Hellstr{\"o}m, L., Silvestrov, S.D.},
  \emph{Commuting
Elements in $q$-Deformed Heisenberg Algebras}, World Scientific,
2000, 256 pp.
\bibitem{Jimbo1D} Jimbo, M., \emph{A $q$-difference analogue of
$U(\g)$ and the Yang--Baxter equation}, Lett. Math. Phys. \textbf{10}
(1985), 63--69.
\bibitem{KodairaSpencerD} Kodaira, K., Spencer, D.C., \emph{On
Deformations of Complex Analytic Structures I-II}, Ann. of Math.
\textbf{67} No 2-3, (1958), 328--466.
\bibitem{KulResD} Kulish, P.P., Reshetikhin, N.Y., \emph{Quantum Linear
    Problem for the Sine-Gordon Equation and Higher Representations},
  Zap. Nauchn. Sem. Leninggrad. Otdel. Mat. Inst. Steklov 101 (1981)
  101--110, Engl. transl. in J. Soviet Math \textbf{23} (1983) no 4.
\bibitem{LarsSilv1D} Larsson, D., Silvestrov, S.D.,
\emph{Quasi-hom-Lie algebras, Central Extensions and
2-cocycle-like identities}, 
  J. Algebra \textbf{288} (2005), 321--344.
\bibitem{LeBruynD} Le Bruyn, L., \emph{Conformal $sl_2$ Enveloping
    Algebras}, Comm. Alg. \textbf{23} no. 4 (1995), 1325--1362.
\bibitem{LeBruynSmithD} Le Bruyn, L., Smith, S.P., \emph{Homogenized
    $sl_2$}, Proc. AMS 118 (1993), 725--730.
\bibitem{LeBruynSmithvdBerghD} Le Bruyn, L., Smith, S.P., Van den Bergh,
  M., \emph{Central Extensions of Three-Dimensional Artin--Schelter
    Regular Algebras}, Math. Zeit. \textbf{222} no. 2 (1996), 171--212.
\bibitem{LeBruynvdBerghD} Le Bruyn, L., Van den Bergh, M., \emph{On Quantum
    Spaces of Lie Algebras}, Preprint, UIA, available @
    \texttt{http://www.math.ua.ac.be/\~{}lebruyn/}.
\bibitem{LevasseurD} Levasseur, T., \emph{Some Properties of
    Non-commutative Regular Graded Rings}, Glasgow Math. J. \textbf{34}
  (1992), 277--300.
\bibitem{LiOys} Huishi, L., Van Oystaeyen, F., \emph{Global
dimension and Auslander regularity of Rees rings}, Bull. Soc.
Math. Belg. \textbf{43} (1991), 59--87.
\bibitem{PiontSilvD} Piontkovski, D., Silvestrov, S.D.,
\emph{Cohomology of 3-dimensional color Lie algebras}, preprint
\texttt{math.KT/0508573}.
\bibitem{Serre1D} Serre, J.P., \emph{Complex semisimple Lie
Algebras}, Springer-Verlag, 2001, 74 pp.
\end{thebibliography}
\end{document}